\documentclass{trans2-l}
\usepackage{amsmath,amsthm,amsfonts,amssymb,url,hyperref}
\usepackage{graphics}
\newcommand{\cyr}{\fontencoding{OT2}\selectfont}
\input{ot2enc.def}
\title{Victor Borisovich Lidskii (1924--2008)}
\date{}
\begin{document}
\maketitle
\begin{center}
\includegraphics{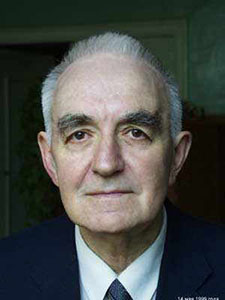}
\end{center}

\

\

\

This volume is dedicated to Victor Borisovich Lidskii who died on 29
July 2008. It is a collection of papers in subject areas related to
Lidskii's work, some of which are written by people who knew him
well. The editors of this volume, Michael Levitin and Dmitri
Vassiliev, are former students of Lidskii, both at undergraduate and
PhD level.

Lidskii was born in Odessa in the Soviet Union on 4 May
1924\footnote{ The correct date of his birth was established only in
the mid 1970s. Lidskii's original birth certificate was lost and
until he got an authorised copy he considered 5 May to be his
birthday and put 5 May in all forms and documents}. As with most men
of his generation, his life was severely affected by the Second
World War. Lidskii finished secondary school on 20 June 1941, two
days before Germany invaded the Soviet Union. After his parents'
divorce he remained with his father and stepmother, so at that
time Lidskii was living in the city of Bobruisk (now Minsk region,
Belarus). His mother lived separately and died in the Minsk ghetto
in 1942.

Lidskii escaped the advancing German troops and found himself, on
his own, in the city of Saratov. Shortly afterwards he was drafted
into the Red Army and sent to a military intelligence school. The
main reason for him being assigned to military intelligence was the
fact that he knew German. Lidskii studied German at school and he
knew some Yiddish as well (which was helpful in learning German).
After graduating as an officer, Lidskii spent the last two years of
the war on the front line serving in field reconnaissance and going
behind enemy lines. Here the German he mastered at intelligence
school proved handy in allowing him, on one occasion, to be taken
for a German by Nazi soldiers when he lead his group back from
behind enemy lines. During these two years he was decorated on four
occasions, see \cite{wikipedia} for details. Lidskii's war ended on
VE Day in Prague and he was demobilised in 1946.

As a decorated officer and member of the Communist Party (he joined
the Party at the front line), Lidskii was offered a job in the NKVD
-- the secret police. Instead, he chose to study mathematics and
enroled as an undergraduate at the Faculty of Mechanics and
Mathematics of Moscow State University. As for many war veterans
studying maths was not easy for him but in the end he proved to be
good enough to be recommended for postgraduate study under the
supervision of I.M.~Gelfand. However, by that time (1951) problems
already started for students of Jewish origin and he was refused a
full-time place at the Graduate School of Moscow State University.
Lidskii managed to secure a teaching post at an engineering college
and did his PhD at Moscow State as an external part-time student.
Despite the difficulties, Lidskii successfully completed his PhD
thesis {\it Questions of spectral theory for systems of second order
differential equations\/} in three years (1954).

In 1954 Lidskii started teaching at the Moscow Institute of Physics
and Technology universally known as FizTech (we write more about
FizTech below) where he worked for the rest of his life. In 1959 he
was awarded a DSc by the Steklov Mathematical Institute for his
thesis {\it Conditions for the completeness of the system of root
subspaces of non-self-adjoint operators with discrete spectra\/} and
shortly afterwards (1961) became a full professor at FizTech. From
1966 until 1988 he also held a part-time research position at the
Institute for Problems in Mechanics of the USSR Academy of Sciences.

V.B.~Lidskii made significant contributions to various branches of
analysis. In this short note we will not attempt to review all his
work. For a survey see~\cite{75} with additional bibliography in
\cite{USSR1,USSR2}. We will only highlight two research topics.

The first --- with which Lidskii's name is most closely associated
--- is the spectral theory of non-symmetric matrices and
non-self-adjoint operators in a Hilbert space. The famous result of
1959 which is widely referred to as Lidskii's Theorem states that
the matrix trace of a trace class operator in a separable Hilbert
space equals its spectral trace (sum of its eigenvalues). Despite
its apparent simplicity it is a deep and difficult achievement, see,
for example, B.~Simon's book~\cite{simon}. Later Lidskii devised a
summation method for series over root vectors of an operator (with
discrete spectrum or compact) which is not necessarily close to
self-adjoint and established sufficient conditions for summability.
In some sense Lidskii's works on non-self-adjoint operators together
with the works of I.~Gohberg and M.G.~Krein were far ahead of their
time. The recent resurgence of interest in non-self-adjoint problems
and the deeper understanding of the underlying difficulties (see
some of the papers in this volume) is strongly influenced by
Lidskii's scientific heritage.

The other, rather unexpected subject, to which Lidskii devoted the
second half of his career, is the \emph{rigorous} analysis of
differential operators arising in elasticity and hydroelasticity.
His favourite research object was the operator of shell theory. This
operator describes the dynamics of a thin elastic structure such as,
say, a hull of a ship or a fuselage of an airplane. In mathematical
terms, dealing with the operator of shell theory means looking at a
self-adjoint Agmon--Douglis--Nirenberg elliptic system of partial
differential equations of mixed order acting on a two-dimensional
manifold with or without boundary. This system depends singularly on
a small asymptotic parameter (relative shell thickness) and can be
reduced to a higher order analogue of the scalar stationary
Schr\"odinger equation with a pseudodifferential potential which
feels in a complicated way the geometry of the manifold. In this
highly nontrivial subject area Lidskii and his collaborators
obtained a collection of delicate rigorous results which at that
time required the use of and contributed to the development of
cutting-edge methods of analysis. Regrettably, these achievements
remain largely unknown: the statements of the problems are too
complicated to attract pure mathematicians and the techniques used
too advanced for the applied mathematics community. Nevertheless,
some of the results buried in Lidskii's works on the subject later
became classical through the work of other authors. For example,
Lidskii was probably the first to prove that the essential spectrum
of a self-adjoint Agmon--Douglis--Nirenberg elliptic operator $A$
(say, on a manifold without boundary) coincides with the set of
values of the spectral parameter $\lambda$ for which the principal
symbol of $A-\lambda I$ is not invertible \cite{AslLid}. The fact
that Lidskii's works on elasticity and hydroelasticity are not
widely known in the West is not surprising given that in Soviet
times he was allowed to travel abroad only once, that only trip
being a trip to Romania.

Lidskii's conversion to shell theory was spontaneous and happened at
Gelfand's seminar at Moscow State. One day in the mid 1960s the
speaker was A.L.~Gol'denveizer, a prominent specialist in shell
theory. After the talk Gelfand stood up and said: ``You, people,
forever keep studying the Laplace operator. Is anyone prepared to
handle something more challenging, like shell theory?'' And Lidskii
was the (only) participant of Gelfand's seminar who accepted this
challenge.

One of the defining aspects of Lidskii's career which is not widely
known was his teaching at FizTech. Here we need to explain the
unique status of this university within the Soviet higher education
system. FizTech was created in 1946 (for the first few years it was
a department of Moscow State University before becoming an
independent institution in 1951) as a result of an initiative of a
group of prominent Soviet scientists which included the future Nobel
laureate P.L.~Kapitsa.  The purpose was to provide elite training of
researchers, mostly in physics, who would later contribute to the
Soviet nuclear and ballistic missile programmes. The university was
based just outside Moscow, in the town of Dolgoprudny, on the site
of the former airship factory of Umberto Nobile. This location
served two purposes: it was a first attempt to create a Western
style university campus where all the students and many of the staff
lived but it also stopped foreigners from visiting as most needed
special permission to leave Moscow. Despite the secrecy and close
association with the military industrial complex, academic freedom
at FizTech was unparalleled. Students at FizTech received additional
extensive training working directly in leading research
establishments of the USSR Academy of Sciences or industry as part
of their undergraduate degree programmes. This gave academics a lot
of freedom in choosing and building their own courses. An
interesting feature of the FizTech ``system'' was the entrance test
which involved not only very strict examinations but also an
informal interview with world leading scientists whose opinion
ultimately determined who was admitted. Unfortunately, as many other
Soviet achievements, this became distorted in later years and strong
anti-Semitic tendencies developed.

The unique features of FizTech attracted not only the best students
from the whole of the USSR but also the best lecturers. The names of
physicists, mathematicians, engineers, chemists and later also
biologists who taught there make an illustrious list. It suffices to
say that Nobel laureates V.~Ginzburg, P.L.~Kapitsa, L.D.~Landau,
A.M.~Prokhorov and N.N.~Semyonov lectured at FizTech. But even
amongst these Lidskii stands out as a true FizTech legend. His
unique passionate teaching style, deep involvement with the subject
and students and strong personality  made him one of the most
popular and loved lecturers, even though he was not by any means a
soft touch as an examiner.

To this day, many years after his peak in the 1960s and 70s,
anecdotes about Lidskii's teaching and personality still circulate
on the web. We cite~\cite{sozer} two such stories here.

Once in a lecture Lidskii gave an example of an everywhere
discontinuous function and sketched its graph. Later on, in the same
lecture, he mentioned that in physics all functions are continuous
and differentiable. He was immediately asked how this agrees with
the quantisation of energy. Lidskii replied that quantisation occurs
precisely because the wave function and its derivative are assumed
to be continuous. At this point he noticed the graph still displayed
on the board. ``This is not a function,'' he said, and stopped and
waved his hand in the air trying to find a proper word, ``but some
sort of pathology!''. He got so angry with this ``pathology'' that
he smashed the board with his fist so that chalk dust flew
everywhere.

On another occasion, in an oral examination for his course, Lidskii
became so upset with a student's answers that in order to calm down
he dashed towards the grand piano standing in the corner of the exam
hall, played angrily some classical piece and only then rejoined the
unfortunate examinee.

It must be said that overall classical music played an important
role in Lidskii's life. It is probably no accident that his youngest
son Mikhail Lidsky, from his second marriage to literary translator
Inna Bernshtein (born 1929), became a well known concert pianist.
Lidskii's first marriage was to Militsa Neuhaus (1929--2008), a
fellow student at Moscow State University and daughter of the famous
pianist Heinrich Neuhaus.

Lidskii also coauthored the legendary exercise book {\it Problems in
Elementary Mathematics} (authors V. Lidsky, L. Ovsyannikov, A.
Tulaikov and M. Shabunin, Mir Publishers, Moscow, 1973) used by
Soviet high school students preparing for university entrance
examinations. Everyone who went on to study maths in the Soviet
Union at university level in the 1960s and 70s associates Lidskii's
name with this exercise book.

The other place where Lidskii worked, part-time, maintaining his
full-time job at FizTech, was the Institute for Problems in
Mechanics (IPMech). IPMech was founded in 1965 by A.Yu.~Ishlinsky, a
distinguished applied mathematician who gained prominence by
developing gyroscopes for the Soviet space programme. Most of the
staff of IPMech had very applied research interests and many were
experimentalists. However, being a well educated and sophisticated
person, Ishlinsky recognised that his institute needs a small group
of ``proper'' mathematicians. The group of mostly part-timers which
came into being as a result included, apart from Lidskii,
O.A.~Oleinik and M.I.~Vishik. IPMech was a purely research
establishment so staff were not obliged to teach and, moreover, the
institute had no undergraduates of its own. Nevertheless, like many
other institutes of the USSR Academy of Sciences, IPMech had close
links with Fiztech: full-timers from IPMech had the opportunity of
teaching part-time at FizTech and full-timers from FizTech had the
opportunity of doing research part-time at IPMech. Some lecture
courses for FizTech students were taught at IPMech and many student
projects and theses were supervised there. Thus, Lidskii's part-time
research job at IPMech was a natural extension of his full-time
teaching job at FizTech.

For several decades Lidskii served as Deputy Editor of
\emph{Functional Analysis and its Applications}, an influential
Soviet journal founded by I.M.~Gelfand. His efforts contributed to
this journal becoming one of the best know in the West and helped to
maintain its reputation in the difficult times following the
collapse of the Soviet Union.

Lidskii had around 20 PhD students most of whom became professional
mathematicians. We were amongst his very last students and remember
fondly not only his academic guidance but his openness and
friendship. Coming to his family apartment for discussion was always
a pleasure. Lunch or dinner was immediately served, questions
regarding family members and events of our lives kindly discussed
and immediate help offered if needed. In many ways Victor Borisovich
Lidskii influenced our careers and shaped our lives.

\

\

\

\hfill Michael Levitin and Dmitri Vassiliev

\end{document}